\documentclass[11pt,a4paper,reqno]{amsart}

\usepackage{latexsym}
\usepackage{amsmath}
\usepackage{amsthm}
\usepackage{amsfonts}
\usepackage{amssymb}
\usepackage{mathrsfs}
\usepackage{graphicx}
\usepackage{float}
\usepackage{booktabs}
\usepackage{array}
\usepackage{xcolor}
\usepackage{tikz}
\usetikzlibrary{arrows,positioning,calc}
\usepackage[colorlinks=true,linkcolor=blue,citecolor=blue,urlcolor=blue]{hyperref}

\allowdisplaybreaks[4]
\newcommand{\Acal}{\mathcal A}

\newcommand{\Ccal}{\mathcal C}
\newcommand{\Dcal}{\mathcal D}

\newcommand{\Ical}{\mathcal I}

\newcommand{\Pcal}{\mathcal P}

\newcommand{\Scal}{\mathcal S}
\newcommand{\Tcal}{\mathcal T}
\newcommand{\Xcal}{\mathcal X}

\newcommand{\proj}{\operatorname{proj}}

\newcommand{\Cone}{\operatorname{Cone}}
\newcommand{\Hom}{\operatorname{Hom}}

\theoremstyle{plain}
\newtheorem{theorem}{Theorem}[section]
\newtheorem{proposition}[theorem]{Proposition}
\newtheorem{lemma}[theorem]{Lemma}
\newtheorem{corollary}[theorem]{Corollary}

\theoremstyle{definition}

\theoremstyle{remark}

\title[Object ideals and idempotent completeness]{Object ideals and idempotent completeness}
\author{Jing He and Yixia Zhang}

\address{School of Mathematics and Statistics, Hunan University of Technology and Business, 410205 Changsha, Hunan, P. R. China}
\email{jinghe1003@163.com}

\address{School of Mathematics and Statistics, Changsha University of Science and Technology, 410114 Changsha, Hunan, P. R. China}
\email{yxzhangmath@163.com}

\date{}
\hypersetup{pdftitle={Object ideals and idempotent completeness},pdfauthor={Panyue Zhou}}

\newcommand{\Fcal}{\mathcal F}
\newcommand{\Jcal}{\mathcal J}
\newcommand{\Ycal}{\mathcal Y}
\newcommand{\Ucal}{\mathcal U}
\newcommand{\Vcal}{\mathcal V}

\newcommand{\Ext}{\operatorname{Ext}}
\newcommand{\End}{\operatorname{End}}
\newcommand{\Ob}{\operatorname{Ob}}
\newcommand{\Ch}{\operatorname{Ch}}
\newcommand{\Kcat}{\mathrm K}
\newcommand{\Gh}{\operatorname{Gh}}
\newcommand{\CoGh}{\operatorname{CoGh}}

\begin{document}

\numberwithin{equation}{section}

\begin{abstract}
We construct an idempotent complete Frobenius exact category with an idempotent complete stable category in which a complete ideal cotorsion pair of object ideals induces a cotorsion pair of objects that fails to be complete. This gives a negative answer to the question of Fu, Guil Asensio, Herzog and Torrecillas under idempotent completeness. The construction uses a lemma on dense triangulated subcategories and an example of Bondarko and Sosnilo in which a weight structure fails to extend to the idempotent completion. The induced cotorsion pair is hereditary, and every object admits an ordinary precover from its left class and an ordinary preenvelope into its right class. We characterize which objects admit special precovers and which admit special preenvelopes.
\end{abstract}

\makeatletter
\@namedef{subjclassname@2020}{\textup{2020} Mathematics Subject Classification}
\makeatother

\maketitle

\section{Introduction}

Fu, Guil Asensio, Herzog and Torrecillas \cite{FGHT13} developed ideal approximation theory in exact categories. An ideal in an additive category consists of additive subgroups of its morphism groups closed under composition on both sides by arbitrary morphisms. For an additive class $\Xcal$, let $[\Xcal]$ denote the ideal of morphisms factoring through an object of $\Xcal$. For an ideal $\Ical$, put
$$
\Ob\Ical=\{X\mid 1_X\in\Ical\}.
$$
An ideal is an object ideal if $\Ical=[\Ob\Ical]$. Fu et al. \cite[Theorem~28 and Question~29]{FGHT13} show that a complete cotorsion pair of objects $(\Xcal,\Ycal)$ induces a complete ideal cotorsion pair $([\Xcal],[\Ycal])$ and ask whether the converse holds.

Ideal approximation has also been developed in extriangulated categories. Zhao and Huang \cite{ZH19} study phantom ideals and introduce object-special ideal precovers. Xu, Fu, Gao and Sun \cite[Theorems~3.5, 4.3 and~5.4]{XFGS26} establish versions of Salce's, Christensen's and Wakamatsu's lemmas in this setting. The existence of object-special ideal precovers does not by definition give special precovers by the associated object class. We explain this distinction for our example in Section \ref{sec:comparison}.

Positive results require additional hypotheses. Sun, Wang and Zhu \cite[Theorem~3.6(2)]{SWZ24} obtain completeness of the object pair when the exact category has enough projectives and injectives and the right ideal is enveloping. Sun, Tan, Wang and Zhu \cite[Corollary~1.3]{STWZ25} prove the implication for Frobenius exact categories whose stable categories are Krull--Schmidt. Zhang and Zhou \cite[Theorem~3.11]{ZZ26} prove the corresponding result for Krull--Schmidt extriangulated categories.

Wang, Wang and Zhu \cite{WWZ26} give a counterexample in an exact category which is not weakly idempotent complete. Ren and Wang \cite[Theorem~1.2 and Question~4.10]{RW26} construct one in a weakly idempotent complete Frobenius exact category with finite dimensional morphism spaces and ask whether the implication holds for idempotent complete exact categories. Their category and its stable category omit direct summands.

We answer this question negatively. All modules are right modules. For an additive category $\Acal$, let $\Ch(\Acal)$ denote the category of unbounded cochain complexes and chain maps, and let $\Kcat(\Acal)$ denote its homotopy category. The morphisms in $\Kcat(\Acal)$ are homotopy classes of chain maps. We use $(C[1])^n=C^{n+1}$ with differential $d_{C[1]}^n=-d_C^{n+1}$. For a ring $S$, let $\proj S$ denote the category of finitely generated projective $S$-modules. Fix
\begin{equation}\label{eq:ring}
R=\mathbb Q[x,y]/(y^2-x^3-x^2),
\hspace{2mm}
\Dcal=\Kcat(\proj R),
\hspace{2mm}
\Tcal=\widetilde{\Dcal}.
\end{equation}
Here $\widetilde{\Dcal}$ is the idempotent completion described by Balmer and Schlichting \cite[Section~1]{BS01}. Its objects are pairs $(B,e)$ with $B\in\Dcal$ and $e^2=e\in\End_{\Dcal}(B)$, and a morphism $(B,e)\to(B',e')$ is a morphism $u:B\to B'$ satisfying $u=e'ue$. The canonical embedding sends $B$ to $(B,1_B)$ and is understood up to isomorphism. Let $\Xcal$ and $\Ycal$ consist of objects isomorphic to complexes in $\Dcal$ concentrated in degrees $\geq0$ and $\leq0$, respectively. Let $\Pcal$ be the category of countably generated projective $R$-modules. Realize $\Tcal$ as the closure of $\Dcal$ under retracts in $\Kcat(\Pcal)$, as justified in Section \ref{sec:exact}. Let $Q:\Ch(\Pcal)\to\Kcat(\Pcal)$ be the canonical quotient functor. Put
$$
\Fcal=\{C\in\Ch(\Pcal)\mid Q(C)\in\Tcal\}.
$$
The exact structure on $\Fcal$ consists of the sequences which split in every degree. The notation $Q^{-1}(\Xcal)$ means the class of complexes whose images lie in $\Xcal$, and similarly for other object classes. Our main result is the following. 

\begin{theorem}{\rm (see Theorems \ref{thm:triangulated} and \ref{thm:exact} for details)}\label{thm:main}
The category $\Fcal$ is an idempotent complete Frobenius exact category and $\underline{\Fcal}\simeq\Tcal$ is idempotent complete. Put
$$
\Ucal=Q^{-1}(\Xcal),\hspace{2mm}
\Vcal=Q^{-1}(\Ycal),\hspace{2mm}
\Ical=[\Ucal],\hspace{2mm}
\Jcal=[\Vcal].
$$
Then $(\Ical,\Jcal)$ is a complete ideal cotorsion pair, with $\Ob\Ical=\Ucal$ and $\Ob\Jcal=\Vcal$. The object cotorsion pair $(\Ucal,\Vcal)$ is hereditary and is not complete. For $C\in\Fcal$, the following conditions are equivalent.
\begin{enumerate}
\item $C$ has a special $\Ucal$-precover.
\item $C$ has a special $\Vcal$-preenvelope.
\item $Q(C)$ belongs to $\Dcal$.
\end{enumerate}
Moreover, $\Ucal$ is precovering and $\Vcal$ is preenveloping in $\Fcal$.
\end{theorem}

The ring and the nonextendible weight structure underlying the construction are due to Bondarko and Sosnilo \cite[Section~3.1]{BS18}. Lemma \ref{lem:dense} connects their example with ideal approximation. It shows that the two halves still generate a complete ideal cotorsion pair after idempotent completion, whereas special object approximations exist precisely on the original triangulated subcategory. We then use the established correspondence between Frobenius categories and their stable categories to obtain the exact counterexample. The contribution is the negative answer under idempotent completeness and the precise description of its failure of special object approximations.
Theorem \ref{thm:main} also applies to the question in extriangulated categories, since exact categories are examples of the framework of Nakaoka and Palu \cite{NP19}. 

The paper is organized as follows. In Section~\ref{sec:triangulated}, we recall the terminology and prove a lemma for dense embeddings. We also construct the triangulated counterexample, including the ring version in Corollary~\ref{cor:family}. In Section~\ref{sec:exact}, we give the Frobenius exact realization. In Section~\ref{sec:comparison}, we compare it with the existing examples and positive results.

\section{Object ideals in the idempotent completion}\label{sec:triangulated}

We first give a construction for a dense triangulated embedding. We then apply it to the homotopy category in \eqref{eq:ring}.

All subcategories are full and closed under isomorphisms. An additive class contains zero and is closed under finite direct sums. A retract of $B$ is an object $A$ with maps $s:A\to B$ and $r:B\to A$ satisfying $rs=1_A$. An additive category is idempotent complete if every idempotent $e:B\to B$ has such a splitting with $sr=e$. It is weakly idempotent complete if every split monomorphism has a cokernel, equivalently every split epimorphism has a kernel. These conventions follow B\"uhler \cite[Sections~6--7]{Buh10}. We call a full triangulated embedding $\Dcal\subseteq\Tcal$ dense if every object of $\Tcal$ is a retract of an object of $\Dcal$.

We use the exact category terminology of B\"uhler \cite[Sections~2 and~13.4]{Buh10}. A conflation is a kernel-cokernel pair belonging to the chosen exact structure. Its first map is an inflation and its second map a deflation. The group $\Ext^1(A,B)$ consists of extension classes of conflations $B\to E\to A$. An exact category $\Ccal$ is Frobenius if it has enough projective and injective objects and these classes coincide. Its stable category $\underline{\Ccal}$ is the quotient by maps factoring through these objects. The suspension $A[1]$ is represented by the cokernel in a conflation $A\to I\to A[1]$ with $I$ injective. This cokernel is called a cosyzygy of $A$.

We recall the notation for ideal orthogonality. In an exact category, take $\mathbb E=\Ext^1$. In a triangulated category, take $\mathbb E(A,B)=\Hom(A,B[1])$. For $i:A\to B$ and $j:C\to D$, write $i\perp j$ when
$$
\mathbb E(i,j):\mathbb E(B,C)\longrightarrow\mathbb E(A,D)
$$
is zero. Thus $\Ical^\perp$ consists of the maps orthogonal to every $i\in\Ical$, and ${}^\perp\Jcal$ is defined dually. An ideal cotorsion pair satisfies $\Ical={}^\perp\Jcal$ and $\Jcal=\Ical^\perp$. We write $j_*\eta=\mathbb E(1,j)(\eta)$ and $i^*\eta=\mathbb E(i,1)(\eta)$ for pushforward and pullback of extensions.

An $\Ical$-precover of $A$ is a morphism $p:P\to A$ in $\Ical$ through which every morphism in $\Ical$ ending at $A$ factors. A $\Jcal$-preenvelope is defined by reversing these arrows. A preenvelope $q:A\to J_A$ is an envelope if every endomorphism $u$ of $J_A$ with $uq=q$ is invertible. An ideal is precovering, preenveloping or enveloping if every object has the corresponding approximation. An $\Ical$-precover occurring in an extension $\delta\in\mathbb E(A,K)$ is special if $\delta=j_*\eta$ for a morphism $j:Z\to K$ in $\Ical^\perp$ and an extension $\eta\in\mathbb E(A,Z)$. A $\Jcal$-preenvelope occurring in an extension $\delta\in\mathbb E(C,A)$ is special if $\delta=i^*\eta$ for a morphism $i:C\to Z$ in ${}^\perp\Jcal$ and an extension $\eta\in\mathbb E(Z,A)$. An ideal cotorsion pair is complete if both kinds of special approximations exist for every object. These are the definitions of Fu et al. \cite{FGHT13} and Breaz and Modoi \cite[Sections~3.2--3.4]{BM21}.

For an object class $\Xcal$, an $\Xcal$-precover of $A$ is a map $p:X_A\to A$ with $X_A\in\Xcal$ through which every map $X\to A$, $X\in\Xcal$, factors. Object preenvelopes are defined dually. A class is precovering or preenveloping if these approximations exist for all objects. We call these ordinary object approximations to distinguish them from special ones. For object classes, orthogonals mean vanishing of $\mathbb E$ on every pair of objects. A cotorsion pair of objects satisfies $\Xcal={}^\perp\Ycal$ and $\Ycal=\Xcal^\perp$. It is complete if every object $M$ occurs in extensions
$$
Y_M\longrightarrow X_M\longrightarrow M,
\hspace{2mm}
M\longrightarrow Y^M\longrightarrow X^M
$$
with all $X$-terms in the left class and all $Y$-terms in the right class. The corresponding maps are special object approximations. In an exact category the pair is hereditary if $\Ext^n(X,Y)=0$ for all $n\geq1$, $X\in\Xcal$ and $Y\in\Ycal$, where $\Ext^n$ denotes Yoneda extensions. In a triangulated category, the displayed extensions are distinguished triangles with their connecting maps. Following Iyama and Yoshino \cite[Definition~2.2]{IY08}, write $\Xcal*\Ycal[1]$ for the class of objects $M$ admitting a triangle $X\to M\to Y[1]\to X[1]$ with $X\in\Xcal$ and $Y\in\Ycal$.

\begin{lemma}\label{lem:dense}
Let $\Dcal\subseteq\Tcal$ be a dense full triangulated subcategory, and $\Xcal,\Ycal\subseteq\Dcal$ be additive subcategories closed under retracts in $\Tcal$. Suppose that
$$
\Hom_{\Tcal}(\Xcal,\Ycal[1])=0
$$
and that every $B\in\Dcal$ admits a triangle
\begin{equation}\label{eq:decomposition}
X_B\xrightarrow{a_B}B\xrightarrow{b_B}Y_B[1]\longrightarrow X_B[1],
\hspace{2mm}X_B\in\Xcal,\ Y_B\in\Ycal.
\end{equation}
Then $([\Xcal],[\Ycal])$ is a complete ideal cotorsion pair in $\Tcal$, with associated object cotorsion pair $(\Xcal,\Ycal)$. The class $\Xcal$ is precovering and $\Ycal$ is preenveloping. For $M\in\Tcal$, a special $\Xcal$-precover exists if and only if $M\in\Dcal$, and a special $\Ycal$-preenvelope exists if and only if $M\in\Dcal$. In particular, $\Xcal*\Ycal[1]=\Dcal$.
\end{lemma}

\begin{proof}
For a class $\Scal$, let $\Gh_{\Scal}$ be the ideal of maps $g:A\to C$ with $gf=0$ for every map $f:S\to A$, $S\in\Scal$. Let $\CoGh_{\Scal}$ be the ideal of maps $f:C\to A$ with $hf=0$ for every map $h:A\to S$, $S\in\Scal$. We first prove
\begin{equation}\label{eq:annihilators}
\Gh_{\Xcal}=[\Ycal[1]],
\hspace{2mm}
\CoGh_{\Ycal[1]}=[\Xcal].
\end{equation}
Choose $A\xrightarrow{s}B\xrightarrow{r}A$ with $B\in\Dcal$ and $rs=1_A$, and take \eqref{eq:decomposition}. If $g:A\to C$ belongs to $\Gh_{\Xcal}$, then $gr a_B=0$. Exactness of the Hom sequence gives $t:Y_B[1]\to C$ such that $gr=t b_B$. This gives the commutative square
$$
\begin{tikzpicture}[baseline=(current bounding box.center),>=stealth]
\node (B) at (0,1.2) {$B$};
\node (Y) at (3.3,1.2) {$Y_B[1]$};
\node (A) at (0,0) {$A$};
\node (C) at (3.3,0) {$C$};
\draw[->] (B) -- node[above] {$b_B$} (Y);
\draw[->] (B) -- node[left] {$r$} (A);
\draw[->] (Y) -- node[right] {$t$} (C);
\draw[->] (A) -- node[below] {$g$} (C);
\end{tikzpicture}
$$
Composing with the section $s$ gives $g=t b_Bs$, so $g$ factors through $Y_B[1]$. The reverse inclusion follows from $\Hom(\Xcal,\Ycal[1])=0$.

For the second equality, let $f:C\to A$ belong to $\CoGh_{\Ycal[1]}$. Since $b_Bsf=0$, exactness gives $u:C\to X_B$ with $sf=a_Bu$. Thus $f=r a_Bu$ factors through $X_B$. Again orthogonality gives the reverse inclusion. This proves \eqref{eq:annihilators}.

Put $\Ical=[\Xcal]$ and $\Jcal=[\Ycal]$. For a map $j$, the condition $j\in\Ical^\perp$ is equivalent to $\Hom(X,j[1])=0$ for every $X\in\Xcal$. Necessity follows by testing on $1_X$. Sufficiency follows by factoring each $i\in\Ical$ through an object of $\Xcal$. Hence $j\in\Ical^\perp$ if and only if $j[1]\in\Gh_{\Xcal}=[\Ycal[1]]$, which is equivalent to $j\in\Jcal$. Similarly, testing $i\in{}^\perp\Jcal$ on the identities $1_Y$ gives $\Hom(i,Y[1])=0$ for all $Y\in\Ycal$. Factoring an arbitrary $j\in\Jcal$ through such a $Y$ proves the converse. Thus $i\in{}^\perp\Jcal$ if and only if $i\in\CoGh_{\Ycal[1]}=[\Xcal]=\Ical$.

Retract closure gives $\Ob\Ical=\Xcal$ and $\Ob\Jcal=\Ycal$, since an identity factors through a class precisely when its object is a retract of that class. Testing the ideal orthogonalities on identity morphisms, as in Sun, Wang and Zhu \cite[Lemma~3.5]{SWZ24}, gives
$$
\Xcal={}^\perp_{\mathbb E}\Ycal,
\hspace{2mm}
\Ycal=\Xcal^{\perp_{\mathbb E}}.
$$

The map $p_A=r a_B:X_B\to A$ is an $\Xcal$-precover. For $f:X\to A$ with $X\in\Xcal$, we have $b_Bsf=0$, so $sf=a_Bu$ and $f=p_Au$. Every morphism of $\Ical$ ending at $A$ factors through such an $f$, so $p_A$ is also an $\Ical$-precover.

For the preenvelope, choose a retraction $A[1]\xrightarrow{s'}B'\xrightarrow{r'}A[1]$ with $B'\in\Dcal$ and $r's'=1_{A[1]}$, together with a decomposition
$$
X_{B'}\xrightarrow{a'}B'\xrightarrow{b'}Y_{B'}[1]\longrightarrow X_{B'}[1].
$$
Put $q_A=(b's')[-1]:A\to Y_{B'}$. Given $h:A\to Y_0$ with $Y_0\in\Ycal$, orthogonality gives $h[1]r'a'=0$. Exactness yields $t:Y_{B'}[1]\to Y_0[1]$ with $h[1]r'=tb'$. Hence $h[1]=tb's'$ and $h=t[-1]q_A$. Thus $q_A$ is a $\Ycal$-preenvelope.

Breaz and Modoi \cite[Example~3.4.2 and Theorem~3.4.1]{BM21} show that a precovering ideal in a triangulated category is special and forms a complete ideal cotorsion pair with its right orthogonal. Consequently $(\Ical,\Jcal)$ is complete.

If $M\in\Dcal$, rotating \eqref{eq:decomposition} gives a triangle $Y_M\to X_M\to M\to Y_M[1]$, hence a special $\Xcal$-precover. Applying \eqref{eq:decomposition} to $M[1]$, shifting by $[-1]$ and rotating gives
$$
M\longrightarrow Y_{M[1]}\longrightarrow X_{M[1]}\longrightarrow M[1],
$$
hence a special $\Ycal$-preenvelope. Conversely, either special approximation gives a triangle whose other two terms lie in $\Dcal$. Closure under shifts and cones therefore places $M$ in $\Dcal$. The same argument shows $\Xcal*\Ycal[1]\subseteq\Dcal$, and \eqref{eq:decomposition} gives the reverse inclusion.
\end{proof}

The distinction in Lemma \ref{lem:dense} can also be read from the connecting map. Complete $p_A:X_B\to A$ to a triangle
$$
K_A\longrightarrow X_B\xrightarrow{p_A}A\xrightarrow{\delta_A}K_A[1].
$$
Since $p_A$ is an $\Xcal$-precover, exactness gives $\Hom(X,\delta_A)=0$ for every $X\in\Xcal$. Hence \eqref{eq:annihilators} gives $\delta_A\in[\Ycal[1]]$. This is a factorization condition on the connecting map. When $A\notin\Dcal$, the triangle cannot have $K_A\in\Ycal$, by Lemma \ref{lem:dense}.

We now apply the lemma to \eqref{eq:ring}. The triangulated structure on $\Tcal=\widetilde{\Dcal}$ and the exactness of the embedding follow from Balmer and Schlichting \cite[Theorem~1.5]{BS01}. Bondarko and Sosnilo \cite[Section~3.1]{BS18} show that $\Dcal$ is not idempotent complete. Schn\"urer \cite[Theorem~3.1]{Schn11} proves that the two halves $\Xcal$ and $\Ycal$ are idempotent complete. Bondarko and Sosnilo \cite[Proposition~1.1.2(1)]{BS18} also record this fact. We use these results without repeating their proofs.

We check the hypotheses involving the ambient category $\Tcal$. An idempotent on an object of either half is a morphism in $\Dcal$, since the embedding is fully faithful. It splits within that half by the cited theorem. Uniqueness of a splitting up to isomorphism then shows that both halves are closed under retracts in $\Tcal$. For $n\geq1$, a representative of an object of $\Ycal[n]$ is concentrated in degrees $\leq-n$. Every component of a chain map from a complex concentrated in degrees $\geq0$ to that representative is zero. Full faithfulness therefore gives
$$
\Hom_{\Tcal}(\Xcal,\Ycal[n])=0\hspace{2mm}(n\geq1).
$$
For a complex $B$, its brutal truncation $\sigma_{\geq m}B$ agrees with $B$ in degrees $\geq m$ and is zero elsewhere, with the remaining differentials inherited from $B$. Define $\sigma_{\leq m}B$ similarly in degrees $\leq m$. Every $B\in\Dcal$ then has the degreewise split sequence
$$
0\longrightarrow\sigma_{\geq0}B\longrightarrow B
\longrightarrow\sigma_{\leq-1}B\longrightarrow0.
$$
It yields a triangle of the form \eqref{eq:decomposition}, with $X_B=\sigma_{\geq0}B$ and $Y_B=(\sigma_{\leq-1}B)[-1]$.

These are the decompositions of the standard weight structure. With the convention of Bondarko and Sosnilo \cite[Definition~1.2.1 and Remark~1.2.3(1)]{BS18}, a weight structure on $\Dcal$ is a pair $(\mathcal L,\mathcal M)$ of additive subcategories closed under retracts, satisfying $\mathcal L[-1]\subseteq\mathcal L$, $\mathcal M[1]\subseteq\mathcal M$, $\Hom(\mathcal L,\mathcal M[1])=0$ and $\Dcal=\mathcal L*\mathcal M[1]$. Here the halves are $(\mathcal L,\mathcal M)=(\Xcal,\Ycal)$.

\begin{theorem}\label{thm:triangulated}
The idempotent complete triangulated category $\Tcal$ in \eqref{eq:ring} admits a complete ideal cotorsion pair $([\Xcal],[\Ycal])$ of object ideals. Its associated object cotorsion pair is not complete. An object $M\in\Tcal$ has a special $\Xcal$-precover, or a special $\Ycal$-preenvelope, exactly when $M\in\Dcal$. Both ordinary approximations exist for every $M\in\Tcal$.
\end{theorem}

\begin{proof}
The preceding observations verify every hypothesis of Lemma \ref{lem:dense}. Since $\Dcal$ is not idempotent complete, its essential image is a proper subcategory of $\Tcal$. Any object in $\Tcal\setminus\Dcal$ therefore has neither of the special object approximations, by that lemma.
\end{proof}

To identify an object without special approximations, choose an idempotent $e:B\to B$ in $\Dcal$ which does not split there. Such an idempotent exists by Bondarko and Sosnilo \cite[Section~3.1]{BS18}. Its image $M=(B,e)$ exists in $\widetilde{\Dcal}$. If $M$ were isomorphic to an object of $\Dcal$, full faithfulness would give a splitting of $e$ in $\Dcal$. Thus $M\notin\Dcal$, and Theorem \ref{thm:triangulated} shows that it has neither special object approximation.

\begin{corollary}\label{cor:family}
Let $S$ be a ring such that $\Kcat(\proj S)$ is not idempotent complete. Its idempotent completion admits a complete ideal cotorsion pair generated by the nonnegative and nonpositive halves. Both ideals are object ideals, and their object pair is not complete.
\end{corollary}

\begin{proof}
Schn\"urer \cite[Theorem~3.1]{Schn11} proves that the two halves are idempotent complete. As above, they are closed under retracts in the idempotent completion, their required Hom groups vanish by support, and brutal truncation supplies the decomposition triangles. Lemma \ref{lem:dense} applies. Since $\Kcat(\proj S)$ is a proper subcategory of its idempotent completion, its object pair is not complete.
\end{proof}

\section{The Frobenius exact counterexample}\label{sec:exact}

We now construct an exact category realizing Theorem \ref{thm:triangulated}. Recall that $\Pcal$ is the category of countably generated projective $R$-modules. This category is idempotent complete and has countable coproducts. Consequently $\Kcat(\Pcal)$ has countable coproducts and is idempotent complete by Neeman \cite[Proposition~1.6.8]{Nee01}. The fully faithful embedding $\Dcal\to\Kcat(\Pcal)$ extends to a fully faithful triangle functor from $\widetilde{\Dcal}$. We identify its image with $\Tcal$. It is a thick subcategory of $\Kcat(\Pcal)$, meaning a full triangulated subcategory closed under retracts.

Let $\Ch(\Pcal)$ carry the exact structure of sequences which split in each degree. A complex is contractible if its identity is chain homotopic to zero. We write $\Cone(f)$ for the mapping cone of a chain map $f:A\to B$, with $\Cone(f)^n=B^n\oplus A^{n+1}$. B\"uhler \cite[Example~13.8]{Buh10} gives the standard Frobenius structure on $\Ch(\Pcal)$, whose stable category is $\Kcat(\Pcal)$ and whose projective and injective objects are the contractible complexes. We restrict this structure to
\begin{equation}\label{eq:exact-category}
\Fcal=\{C\in\Ch(\Pcal)\mid Q(C)\in\Tcal\}.
\end{equation}
Equivalently, its objects are the complexes which are retracts, in the homotopy category, of complexes of finitely generated projective modules.

\begin{proposition}\label{prop:model}
With the induced degreewise split exact structure, $\Fcal$ is an idempotent complete Frobenius exact category. Its projective and injective objects are the contractible complexes, and $Q$ induces a triangle equivalence $\underline{\Fcal}\simeq\Tcal$.
\end{proposition}

\begin{proof}
A degreewise split conflation in $\Ch(\Pcal)$ induces a triangle in $\Kcat(\Pcal)$. If its end terms belong to $\Fcal$, their images under $Q$ lie in $\Tcal$. The middle term also has image in $\Tcal$, since $\Tcal$ is triangulated. Thus $\Fcal$ is extension closed and inherits an exact structure by B\"uhler \cite[Lemma~10.20]{Buh10}. Every contractible complex belongs to $\Fcal$, because its image under $Q$ is zero.

An idempotent chain map on $C\in\Fcal$ splits degreewise in $\Pcal$. Its image is a subcomplex $C'$ and $Q(C')$ is a retract of $Q(C)$. Since $\Tcal$ is closed under retracts, $C'$ belongs to $\Fcal$. Thus $\Fcal$ is idempotent complete.

For $C\in\Fcal$, the standard cone sequences are
$$
0\longrightarrow C[-1]\longrightarrow\Cone(1_{C[-1]})\longrightarrow C\longrightarrow0,
$$
$$
0\longrightarrow C\longrightarrow\Cone(1_C)\longrightarrow C[1]\longrightarrow0.
$$
All terms belong to $\Fcal$, because $\Tcal$ is closed under shifts and both cones are contractible. The contractible complexes remain projective and injective under restriction of the exact structure. These sequences therefore give enough projectives and injectives in $\Fcal$. Any projective object is a retract of the contractible middle term of the first sequence, and any injective object is a retract of the middle term of the second. Both kinds of retract are contractible. Hence the projective and injective objects of $\Fcal$ coincide and are exactly the contractible complexes.

Finally, the kernel ideal of $Q$ consists of the maps factoring through contractible complexes by B\"uhler \cite[Remark~9.11]{Buh10}. All these complexes belong to $\Fcal$, so the restriction of $Q$ induces a fully faithful functor $\underline{\Fcal}\to\Kcat(\Pcal)$. Its essential image is $\Tcal$ by \eqref{eq:exact-category}. The inclusion into $\Ch(\Pcal)$ is exact and preserves projective and injective objects. It therefore induces a triangle functor for the standard stable structures described by B\"uhler \cite[Theorem~13.7 and Example~13.8]{Buh10}. This gives the required triangle equivalence.
\end{proof}

Set $\Ucal=Q^{-1}(\Xcal)$ and $\Vcal=Q^{-1}(\Ycal)$. Both are closed under retracts and contain all contractible complexes. For an ideal $\mathcal L$ in $\Tcal$, write $Q^{-1}(\mathcal L)=\{f\mid Q(f)\in\mathcal L\}$. Fullness of $Q$ and the existence of representatives in these classes give $Q([\Ucal])=[\Xcal]$ and $Q([\Vcal])=[\Ycal]$. Since both ideals contain all maps factoring through projective objects, Sun, Tan, Wang and Zhu \cite[Lemma~3.2]{STWZ25} give
\begin{equation}\label{eq:lifted-ideals}
Q^{-1}([\Xcal])=[\Ucal],
\hspace{2mm}
Q^{-1}([\Ycal])=[\Vcal].
\end{equation}

\begin{theorem}\label{thm:exact}
The category $\Fcal$ and the classes $\Ucal,\Vcal$ satisfy all conclusions of Theorem \ref{thm:main}.
\end{theorem}

\begin{proof}
Proposition \ref{prop:model} gives both idempotent completeness assertions. Saor\'in and \v{S}\v{t}ov\'i\v{c}ek \cite[Proposition~3.3 and its proof]{SS11} give the natural Frobenius isomorphism
\begin{equation}\label{eq:ext-stable}
\Ext^1_{\Fcal}(A,B)\cong\Hom_{\Tcal}(Q(A),Q(B)[1]).
\end{equation}
Every object of $\Xcal$ has a representative in $\Ucal$, and every object of $\Ycal$ has a representative in $\Vcal$. Consequently \eqref{eq:ext-stable} and Lemma \ref{lem:dense} identify the object orthogonals as $\Ucal={}^\perp\Vcal$ and $\Vcal=\Ucal^\perp$. To identify the ideal orthogonals, test on the identities of these representatives and use naturality of \eqref{eq:ext-stable}. Together with \eqref{eq:lifted-ideals}, this gives
$$
\begin{aligned}
\relax [\Ucal]^\perp
&=Q^{-1}([\Xcal]^\perp)=Q^{-1}([\Ycal])=[\Vcal],\\
{}^\perp[\Vcal]
&=Q^{-1}({}^\perp[\Ycal])=Q^{-1}([\Xcal])=[\Ucal].
\end{aligned}
$$
Retract closure gives $\Ob[\Ucal]=\Ucal$ and $\Ob[\Vcal]=\Vcal$.

By \eqref{eq:lifted-ideals}, the stable quotient of this ideal cotorsion pair is the complete pair of Theorem \ref{thm:triangulated}. Both ideals contain the identities of all projective objects. The completeness correspondence of Sun, Tan, Wang and Zhu \cite[Corollary~1.2]{STWZ25} therefore shows that $([\Ucal],[\Vcal])$ is complete. The same result implies that $[\Ucal]$ is precovering and $[\Vcal]$ is preenveloping. By Sun, Wang and Zhu \cite[Lemma~3.2(1)]{SWZ24}, it follows that $\Ucal$ is precovering and $\Vcal$ is preenveloping.

For the special object approximations, we use the Frobenius correspondence of Saor\'in and \v{S}\v{t}ov\'i\v{c}ek \cite[Proposition~3.3]{SS11}, in the form for individual objects stated by Ren and Wang \cite[Lemma~4.4]{RW26}. Under the equivalence $\underline{\Fcal}\simeq\Tcal$, it gives
$$
\begin{aligned}
C\text{ has a special }\Ucal\text{-precover}
&\ \Longleftrightarrow\ Q(C)\in\Xcal*\Ycal[1],\\
C\text{ has a special }\Vcal\text{-preenvelope}
&\ \Longleftrightarrow\ Q(C)[1]\in\Xcal*\Ycal[1].
\end{aligned}
$$
The hypotheses of this correspondence hold because $(\Ucal,\Vcal)$ is a cotorsion pair, both classes contain all projective and injective objects, and they are the full inverse images of their stable classes. Lemma \ref{lem:dense} gives $\Xcal*\Ycal[1]=\Dcal$. Since $\Dcal$ is closed under both shifts, each right hand condition is equivalent to $Q(C)\in\Dcal$. This proves the three equivalent conditions in Theorem \ref{thm:main}. Essential surjectivity of $Q$ and the proper inclusion $\Dcal\subsetneq\Tcal$ provide objects for which both special approximations fail, so the object pair is not complete.

Finally, the cone conflations in Proposition \ref{prop:model} identify the stable image of the $m$th cosyzygy of $V$ with $Q(V)[m]$. Dimension shifting through their injective middle terms and then applying \eqref{eq:ext-stable} gives
$$
\Ext^n_{\Fcal}(U,V)\cong\Hom_{\Tcal}(Q(U),Q(V)[n])=0
\hspace{2mm}(n\geq1)
$$
for $U\in\Ucal$ and $V\in\Vcal$. The last vanishing was verified by support in Section \ref{sec:triangulated}. Hence the object cotorsion pair is hereditary.
\end{proof}

\section{Comparison with existing examples}\label{sec:comparison}

Theorem \ref{thm:exact} identifies exactly which objects lack special approximations. They are the representatives of objects in $\Tcal\setminus\Dcal$. We finish by comparing this construction with earlier examples and the hypotheses of the positive results.

Wang, Wang and Zhu \cite{WWZ26} use an extension closed subcategory of a module category defined by an additive numerical invariant. Their category is not weakly idempotent complete. Ren and Wang \cite[Sections~3--4]{RW26} use bounded complexes of finite dimensional vector spaces with even total cohomology dimension. Their category is weakly idempotent complete but misses direct summands, and idempotent completion removes that obstruction. In Theorem \ref{thm:main}, both the exact category and its stable category are idempotent complete. Special approximations fail because $\Xcal*\Ycal[1]=\Dcal$ is a proper subcategory of $\Tcal$, although $\Xcal$ and $\Ycal$ are closed under retracts.

Our morphism spaces need not be finite dimensional over $\mathbb Q$. For example, $R[0]$, the complex with $R$ in degree zero and zero elsewhere, has endomorphism ring $R$, which is infinite dimensional over $\mathbb Q$. Moreover, $R$ is a domain with distinct maximal ideals $(x,y)$ and $(x-3,y-6)$, so it is not local. Thus $R[0]$ is indecomposable and has a nonlocal endomorphism ring. A Krull--Schmidt category is one in which every object is a finite direct sum of objects with local endomorphism rings, as in Krause \cite[Section~4]{Kra15}. Consequently $\Tcal$ is not Krull--Schmidt. The results of Sun, Tan, Wang and Zhu \cite[Corollary~1.3]{STWZ25} and Zhang and Zhou \cite[Theorem~3.11]{ZZ26} therefore do not apply.

There is also a distinction between the two notions of special approximation. Zhao and Huang \cite[Definition~5.1]{ZH19} call an $\Ical$-precover $q:E\to C$ object-special if it occurs in a conflation $K\to E\xrightarrow{q}C$ with $K\in\Ob(\Ical^\perp)$. This definition requires $q\in\Ical$, but does not require $E\in\Ob\Ical$. Zhao and Huang \cite[Proposition~5.2]{ZH19} prove that a special precovering ideal is object-special precovering whenever its right orthogonal is an object ideal. Applied to Theorem \ref{thm:main}, their result gives, for every $C\in\Fcal$, a conflation
$$
0\longrightarrow V_C\longrightarrow E_C\xrightarrow{q_C}C\longrightarrow0,
\hspace{2mm}V_C\in\Vcal,\hspace{2mm}q_C\in[\Ucal].
$$
If $Q(C)\notin\Dcal$, then $E_C\notin\Ucal$. Otherwise this conflation would be a special $\Ucal$-precover, contrary to Theorem \ref{thm:main}. Thus the factorization of $q_C$ through an object of $\Ucal$ does not place its domain in $\Ucal$. Dually, object-special ideal preenvelopes exist for every object, although special $\Vcal$-preenvelopes fail outside $Q^{-1}(\Dcal)$. Our example therefore satisfies the conclusion of Zhao and Huang while separating it from completeness of the associated object cotorsion pair.

For a cotorsion pair in a Frobenius category, cohereditary means that its stable left class is closed under $[1]$. Ren and Wang \cite[Corollary~4.5]{RW26} obtain a positive result for such pairs with idempotent complete stable category. They use the criterion of Saor\'in and \v{S}\v{t}ov\'i\v{c}ek \cite[Proposition~3.11]{SS11} for a precovering full subcategory closed under extensions, retracts and $[1]$. Our left class satisfies $\Xcal[-1]\subseteq\Xcal$. However, $R[0]\in\Xcal$ and $R[0][1]\notin\Xcal$, because the latter has nonzero cohomology in degree $-1$. Thus our hereditary pair is not cohereditary. The right ideal $[\Vcal]$ is also not enveloping, since otherwise the theorem of Sun, Wang and Zhu \cite[Theorem~3.6(2)]{SWZ24} would make the object pair complete.

Finally, we compare with ideal torsion pairs. Following Zhu, Fu, Herzog and Schlegel \cite[Definition~2, Theorem~19 and Example~23]{ZFHS25}, an ideal torsion pair $(\mathcal L,\mathcal M)$ is a pair of ideals that are each other's annihilators for composition. Thus every composable $\ell\in\mathcal L$ and $m\in\mathcal M$ satisfy $m\ell=0$, and each ideal contains all maps with this annihilation property against the other. In a triangulated category the pair is complete if every $A$ occurs in a triangle $P\xrightarrow{\ell}A\xrightarrow{m}N\to P[1]$ with $\ell\in\mathcal L$ and $m\in\mathcal M$.

In the version specified in the bibliography, Zhu, Fu, Herzog and Schlegel \cite[Corollary~56]{ZFHS25} assert completeness of the object pair for every complete ideal cotorsion pair of object ideals in a Frobenius category. Theorem \ref{thm:main} shows that this conclusion requires an additional hypothesis, even when both the exact category and its stable category are idempotent complete. In our stable category we have
$$
([\Xcal],[\Ycal[1]])\text{ is a complete ideal torsion pair in }\Tcal,
$$
$$
\Xcal*\Ycal[1]=\Dcal\subsetneq\Tcal.
$$
Here ideal torsion pairs use Hom-orthogonality, whereas the ideal cotorsion pair in Theorem \ref{thm:triangulated} is $([\Xcal],[\Ycal])$ and uses $\mathbb E=\Hom(-,-[1])$. The first assertion follows from Lemma \ref{lem:dense} and the precover triangles, using Zhu, Fu, Herzog and Schlegel \cite[Example~23]{ZFHS25}. In a triangle $X_B\to A\xrightarrow{\delta_A}K_A[1]\to X_B[1]$, the map $\delta_A$ factors through $\Ycal[1]$, although $K_A[1]$ need not belong to $\Ycal[1]$. An object torsion pair $(\Xcal,\Ycal[1])$ in the sense of Iyama and Yoshino \cite[Definition~2.2]{IY08} requires both $\Hom(\Xcal,\Ycal[1])=0$ and $\Tcal=\Xcal*\Ycal[1]$. Our example satisfies the orthogonality but fails the decomposition condition. Thus the passage from complete object ideals to an object torsion pair in their proof requires an additional hypothesis. The completeness transfer used here is the separate result of Sun, Tan, Wang and Zhu \cite[Corollary~1.2]{STWZ25}.

\vspace{3mm}

\noindent\textbf{Acknowledgements:}
Jing He is supported by the National Natural
Science Foundation of China (Grant No. 12401045).

\vspace{3mm}

\noindent{\bf Statement on AI:} ChatGPT was used only for language editing and improving the clarity of the manuscript. All mathematical results, arguments, and proofs were developed and verified by the authors.

\end{document}